\documentclass[onecolumn]{autart}
\usepackage{}    

\usepackage{graphicx}           
\usepackage{dcolumn}            
\usepackage{bm}                 

\usepackage{graphics}           
\usepackage{epsfig}             
\usepackage{times}              
\usepackage[fleqn]{amsmath}
\usepackage{amssymb}
\usepackage{extarrows}
\usepackage{amsfonts}
\usepackage{mathrsfs}
\usepackage{fancyhdr}
\usepackage{color}
\usepackage{subfigure}
\usepackage{enumerate}
\usepackage [latin1]{inputenc} %
\allowdisplaybreaks[4] 

\newtheorem{theorem}{Theorem}
\newtheorem{lemma}[theorem]{Lemma}

\newtheorem{proposition}[theorem]{Proposition}

\newtheorem{definition}{Definition}

\newtheorem{remark}{Remark}

\begin{document}

\begin{frontmatter}

\title{A Note on the Maximum Principle-based Approach for ISS Analysis of Higher Dimensional Parabolic PDEs with Variable Coefficients}

\author{Jun Zheng$^{1}$}\ead{zhengjun2014@aliyun.com},
\author{Guchuan Zhu$^{2}$}\ead{guchuan.zhu@polymtl.ca},

\address{$^{1}${School of Mathematics, Southwest Jiaotong University,
        Chengdu 611756, Sichuan, China}\\
        $^{2}$Department of Electrical Engineering, Polytechnique Montr\'{e}al, P.O. Box 6079, Station Centre-Ville, Montreal, QC, Canada H3T 1J4}

\begin{abstract}
This paper presents a maximum principle-based approach in the establishment of input-to-state stability (ISS) for a class of nonlinear parabolic partial differential equations (PDEs) over higher dimensional domains with variable coefficients and different types of nonlinear boundary conditions. Technical development on ISS analysis of the considered systems is detailed, and an example of establishing ISS estimates for a nonlinear parabolic equation with, respectively, a nonlinear Robin boundary condition and a nonlinear Dirichlet boundary condition is provided to illustrate the application of the developed method.
\end{abstract}
\end{frontmatter}
 \section{Introduction}\label{Sec: Introduction}
Since the last decade, the ISS theory for infinite dimensional systems governed by partial differential equations (PDEs) has drawn much
attention in the literature of PDE control. {}{A comprehensive survey on this topic is presented in \cite{Mironchenko:2019b}}. It is worth noting that the extension of the notion of ISS for finite dimensional systems  originally introduced by Sontag in the late {}{1980's} to infinite dimensional systems with distributed in-domain disturbances is somehow straightforward, while the investigation on the ISS properties with respect to (w.r.t.) boundary disturbances is much more challenging. In recent years, different methods have been developed for ISS analysis of PDE systems with boundary disturbances, including, e.g.:
\begin{enumerate}[(i)]
\item 
    {}{the semigroup and admissibility methods for ISS of certain linear or nonlinear parabolic PDEs \cite{Jacob:2018CDC,Jacob:2019,jacob2016input,Jacob:2018_SIAM,Jacob:2018_JDE,Schwenninger:2019}};\label{item 1}
   \item the approach of spectral decomposition and finite-difference scheme for ISS of PDEs governed by Sturm-Liouville operators \cite{Karafyllis:2014,Karafyllis:2016,Karafyllis:2016a,karafyllis2017siam,Karafyllis:2017};\label{item 3}
   \item {}{the Riesz-spectral approach for ISS of Riesz-spectral systems \cite{Lhachemi:201901,Lhachemi:201902};\label{item 4}}

               \item the monotonicity-based method for ISS of {certain nonlinear PDEs with Dirichlet boundary disturbances} \cite{Mironchenko:2019};\label{item 5}
   \item the method of De~Giorgi iteration for ISS of certain {nonlinear PDEs with Dirichlet boundary disturbances} \cite{Zheng:201702,Zheng:201803};\label{item 6}
   \item the application of variations of Sobolev embedding inequalities for ISS of certain nonlinear PDEs with Robin or Neumann boundary disturbances \cite{Zheng:201804,Zheng:201803,Zheng:201802};\label{item 7}

   \item the maximum principle-based approach for ISS of certain {nonlinear PDEs with different types of boundary conditions} \cite{Zheng:2019a,Zheng:2019b}.\label{item 8}

\end{enumerate}
{}{Although a rapid progress on ISS theory has been obtained,} it is still a challenging issue for ISS analysis of nonlinear PDE systems defined over higher dimensional domains with variable coefficients and different types of nonlinear boundary conditions. For example, the methods in 
{}{(\ref{item 1})} can be applied to certain linear or nonlienar PDEs, {}{while it may be difficult to apply them to non-diagonal systems as the one given by, e.g., \eqref{LPE1'} and \eqref{Robin} of this paper}. 
The methods in (ii) and (iii) are effective for ISS analysis of linear {}{PDE systems} over one dimensional domains. {}{Whereas, these approaches may involve heavy
computations for nonlinear PDEs or PDEs on multidimensional spatial domains.}
{}{The methods in (iv)-(vi) are suitable for ISS analysis of parabolic PDEs with Dirichlet or Robin or Neumann boundary disturbances, while they cannot be used for PDEs with mixed boundary conditions.}

{}{It has been demonstrated in \cite{Zheng:2019a} and \cite{Zheng:2019b} that the method in (vii) is applicable for ISS analysis of certain nonlinear parabolic PDEs, with different types of boundary disturbances, or over higher dimensional domains.}
Therefore, the aim of this paper is put on the application of this approach to ISS analysis for a class of nonlinear parabolic PDEs defined over higher dimensional domains with variable coefficients under different types of nonlinear boundary conditions simultaneously.

The proposed method for achieving the ISS estimates of the solutions will be based on the the Lyapunov method and the maximum estimates
for nonlinear parabolic PDEs with nonlinear boundary conditions. Specifically, we set up in the first step several maximum estimates of the solutions to the considered nonlinear parabolic PDEs with a nonlinear Robin or Dirichlet boundary condition by means of the weak maximum principle. In the second step, {}{applying the technique of splitting as in \cite{Fabre:1995,Zheng:201702,Zheng:2019a,Zheng:2019b}, we consider a nonlinear equation with the initial data free} and establish the maximum estimate of the solution (denoted by $v$) according to the result obtained in the first step. By denoting the solution of the {}{target} system by $u$, then in the third step, we establish the $L^2$-estimate of $u-v$ by the Lyapunov method. Finally, the ISS estimate of the {}{target} system in $L^2$-norm, i.e., the estimate of $u$ in $L^2$-norm, is guaranteed by the maximum estimate of $v$ and the $L^2$-estimate of $u-v$. {}{It's worthy noting that combining with other approaches or techniques, the Lyapunov method was also applied for the ISS analysis of PDE systems in \cite{Jacob:2019} by constructing non-coercive Lyapunov functions based on ISS characterizations devised in \cite{Mironchenko:2018}, and in \cite{Schwenninger:2019,Tanwani:2017} and the literature mentioned in (iv)-(vii) by constructing coercive Lyapunov functions.}

In the rest of the paper, Section~\ref{Maximum principle} presents the problem statement, the basic assumptions, and the main result. By the weak maximum principle, some maximum estimates for nonlinear parabolic PDEs with nonlinear Robin and Dirichlet boundary conditions are proved respectively in Section~\ref{Sec: Maximum estimate}. ISS analysis of nonlinear parabolic PDEs with different boundary conditions are detailed in Section~\ref{Sec: EISS estimates}. In order to illustrate the application of the approach presented in this paper, an example of ISS estimates for a parabolic equation with respectively a nonlinear Robin boundary condition and a {}{nonlinear} Dirichlet boundary condition is provided in Section~\ref{Sec: example}, followed by some concluding remarks given in Section~\ref{Sec: Conclusion}.

\emph{\textbf{Notations:}} In this paper, $\mathbb{R}_+$ denotes the set of positive real numbers and $\mathbb{R}_{\geq 0} :=  {}{\{0\}}\cup\mathbb{R}_+$.

Let $B_{R}$ be a ball in $\mathbb{R}^n (n\geq 1)$ with the centre at $0$ and a radius $R>0$, i.e.,
$B_R=\{x\in\mathbb{R}^n | |x|<R \}$. We denote by $\partial B_{R}$ and $\overline{B}_{R}$ the boundary and the closure of $B_{R}$, respectively. Denote by $|B_R|$ the $n$-dimensional Lebesgue measure of $B_R$,  i.e., $|B_R|=\frac{\pi^{\frac{n}{2}}R^n}{\frac{n}{2}\Gamma(\frac{n}{2} )}$ with the Gamma function $\Gamma (\cdot)$ defined on $\mathbb{R}$.
Denote by $ |\partial B_R|$ the $(n-1)$-dimensional Lebesgue measure of $\partial B_R$, i.e., $|\partial B_R|=\frac{2\pi^{\frac{n}{2}}}{\Gamma(\frac{n}{2} )}R^{n-1}$.

For any $T>0$, {}{let} $ Q_T=B_{R}\times (0,T)$ and 
$\partial_pQ_T=(\partial B_{R}\times (0,T))\cup (\overline{B}_{R}\times \{{}{0}\})$.

We use $\|\cdot\|$ to denote the norm $ \|\cdot\|_{L^2(B_R)}$ in $L^{2}(B_R)$.

Let $\mathcal {K}=\{\gamma : \mathbb{R}_{\geq 0} \rightarrow \mathbb{R}_{\geq 0}|\ \gamma(0)=0,\gamma$ is continuous, strictly increasing$\}$;
$ \mathcal {K}_{\infty}=\{\theta \in \mathcal {K}|\ \lim\limits_{s\rightarrow\infty}\theta(s)=\infty\}$;
$ \mathcal {L}=\{\gamma : \mathbb{R}_{\geq 0}\rightarrow \mathbb{R}_{\geq 0}|\ \gamma$ is continuous, strictly decreasing, $\lim\limits_{s\rightarrow\infty}\gamma(s)=0\}$;
$ \mathcal {K}\mathcal {L}=\{\mu : \mathbb{R}_{\geq 0}\times \mathbb{R}_{\geq 0}\rightarrow \mathbb{R}_{\geq 0}|\ \mu \in \mathcal {K}, \forall t \in \mathbb{R}_{\geq 0}$, and $\mu (s,\cdot)\in \mathcal {L}, \forall s \in {\mathbb{R}_{+}}\}$.

%

\section{Problem Setting and Main Result}\label{Maximum principle}

{}{Given the following functions:
\begin{subequations}\label{continuity}
\begin{align}
&a,c\in C^2(\overline{B}_R\times \mathbb{R}_{\geq 0};\mathbb{R}_+),\\
&b_i\in C^2(\overline{B}_R\times \mathbb{R}_{\geq 0};\mathbb{R}),i=1,2,...,n,\\
&h\in C^2(\overline{B}_R\times \mathbb{R}_{\geq 0}\times \mathbb{R}; \mathbb{R}),\psi\in C^2(\mathbb{R};\mathbb{R}),\\
&f,d\in C^2(\overline{B}_R\times \mathbb{R}_{\geq 0};\mathbb{R}),\phi\in C^2(\overline{B}_R; \mathbb{R}),
\end{align}
\end{subequations}
we consider the following nonlinear parabolic equation with variable coefficients:
\begin{subequations}\label{LPE1'}
\begin{align}
 L_t[u](x,t)+N[u](x,t)&=f(x,t), (x,t)\in  B_R\times \mathbb{R}_+,\\
\mathscr{B}[u](x,t)&=d(x,t), (x,t)\in \partial B_{R} \times \mathbb{R}_+,\\
u(x,0)&=\phi(x), x\in  B_{R},
\end{align}
\end{subequations}
where $ L_t[u]:=u_t-\text{div} \ (a\nabla u)
+\textbf{b}\cdot \nabla u +cu
$ with $\textbf{b}:=(b_1,b_2,\ldots,b_n)$, $N[u]:=h(\cdot,\cdot,u)$ is the nonlinear term of the equation, $ \mathscr{B}[u]$ is given by:
\begin{align}\label{Robin}
\mathscr{B}[u]:=\frac{\partial u}{\partial\bm{\nu}}+\psi(u),
\end{align}
or
\begin{align}\label{Neumann}
\mathscr{B}[u]:=\psi\bigg(\frac{\partial u}{\partial\bm{\nu}}\bigg),
\end{align}
or
\begin{align}\label{Dirichlet}
\mathscr{B}[u]:=\psi(u),
\end{align}
with $\bm{\nu}=\frac{1}{R}x=\frac{1}{R}(x_1,x_2,\ldots,x_n)$, which is the unit outer normal vector at the point $x\in \partial B_R$.}

{}{In general, $f$ and $d$ represents the distributed in-domain disturbance and boundary disturbance, respectively. \eqref{Robin}, \eqref{Neumann} and \eqref{Dirichlet} represent the nonlinear Robin, Neumann and Dirichlet boundary condition, respectively.}

{}{Throughout this paper, without special statements, we always denote by $x,t$ respectively the first and second variable (if any) of the functions $a,b_i \ (i=1,2,\ldots,n),c,h,f,d,\phi$. Moreover, we always assume that $a,b_i \ (i=1,2,\ldots,n),c,h,\psi,f,d,\phi$ are given by \eqref{continuity} and satisfy for some $\underline{a},\overline{a},\underline{b},\overline{b},\underline{c}\in \mathbb{R}_{+}$:}
\begin{subequations}\label{ABCCBA}
\begin{align}
&0<\underline{a}\leq {}{a\leq \overline{a}\ \text{and}\ |\nabla a|\leq \overline{a}  \ \text{in}\   {B}_R\times \mathbb{R}_{+}},\label{A}\\
&0\leq \underline{b}\leq |\textbf{b}|+|\text{div}\ \textbf{b} |\leq \overline{b} {}{ \ \text{in}\   {B}_R\times \mathbb{R}_{+},}\label{B}\\
&0< \underline{c}\leq c{}{ \ \text{in}\   {B}_R\times \mathbb{R}_{+},}\label{C}\\
&\overline{b}\left(1+2C^2_{\text{Trace}}\right)<2\underline{c},\ \overline{b}C^2_{\text{Trace}}< \underline{a},\label{ABC}
\end{align}
\end{subequations}
 {}{where $C_{\text{Trace}}$ is the best constant of the trace embedding inequality given by the Trace Theorem in the appendix, and
\begin{subequations}\label{Hpsi}
\begin{align}
&N[u](x,t)<N[v](x,t),N[w](x,t)+N[-w](x,t)\geq 0 \ \text{and}\ N[0](x,t)=0,\label{H}\\
&\psi(u)<\psi(v),\psi(w)+\psi(-w)\geq 0\ \text{and}\ \psi(0)=0,\label{psi}
\end{align}
\end{subequations}
for all $(x,t)\in B_{R}\times \mathbb{R}_+$ and all $u,v,w\in \mathbb{R}$ with $u<v$.}

 {}{Furthermore, we impose the following compatibility condition:
\begin{align}\label{compatibility 1}
\mathscr{B}[\phi](x,0)=d(x,0)=0,x\in\partial B_R.
\end{align}}
Let 
$\mathbb{Y} := C^2(\overline{B}_R\times\mathbb{R}_{\geq 0};\mathbb{R})$, $\mathbb{D}_0: =  \{d\in \mathbb{Y}|d(\cdot,0)=0\ \text{on}\ \partial B_R\}$ and $\mathbb{U}_0: =  \{u\in \mathbb{Y}| \mathscr{B}[u](\cdot,0)=d(\cdot,0)\ \text{on}\ \partial B_R \ \text{for}\ d\in \mathbb{D}_0\}$.


\begin{definition}
System~\eqref{LPE1'} is said to be input-to-state stable (ISS) in $L^2$-norm w.r.t. the boundary disturbance $d\in \mathbb{D}_0$, the in-domain disturbance $f\in \mathbb{Y}$ and the states in $ \mathbb{U}_0$, if there exist functions $\beta\in \mathcal {K}\mathcal {L}$ and $ \gamma_0,\gamma_1\in \mathcal {K}$ such that the solution of \eqref{LPE1'} satisfies for any $T>0$:
\begin{align}
    \|u(\cdot,T)\|\leq& \beta(\|\phi\|,T)+\gamma_0 \left(\|d\|_{L^{\infty}(\partial B_R\times (0,T))}\right)
      +\gamma_1\left(\|f\|_{L^{\infty}(Q_T)}\right).\label{Eq: ISS def}
\end{align}
Moreover, System~\eqref{LPE1'} is said to be exponential input-to-state stable (EISS) in $L^2$-norm w.r.t. the boundary disturbance $d\in \mathbb{D}_0$, the in-domain disturbance $f\in \mathbb{Y}$ and the states in $ \mathbb{U}_0$, if $\beta( \|\phi\|,T)$ can be chosen as $ M_0\e^{-\lambda T}\|\phi\|$ with certain constants $M_0,\lambda>0$ in \eqref{Eq: ISS def}.
\end{definition}
The main result {}{of this paper is stated in} the following theorem.
\begin{theorem}\label{main result}
{}{System \eqref{LPE1'} with \eqref{Robin} (or \eqref{Neumann}, or \eqref{Dirichlet})} is EISS w.r.t. the boundary disturbance $d\in \mathbb{D}_0$, the in-domain disturbance $f\in \mathbb{Y}$ and the states in $ \mathbb{U}_0$ having the estimate given in \eqref{2801} (or \eqref{2802}, or \eqref{2803}).
\end{theorem}
\begin{remark}
\begin{enumerate}[(i)]
\item {As pointed out in \cite{Levine:1974}, for a heat conduction problem} the nonlinear boundary conditions can be seen as a nonlinear radiation law prescribed on the boundary of the material body.
\item {}{By \cite[Theorem 6.1 and 7.4, Chapter V]{Ladyzenskaja:1968}, system~\eqref{LPE1'} admits a unique solution ${}{u\in C^{2,1}(\overline{Q}_T)}$ for any $T>0$. Moreover, every system appearing in this paper admits a unique solution belonging to $C^{2,1}(\overline{Q}_T)$.}


 \end{enumerate}
\end{remark}
\begin{remark}
\begin{enumerate}[(i)]
 \item {}{As $\psi$ is  invertible, the nonlinear boundary conditions \eqref{Neumann} and \eqref{Dirichlet} are equivalent to the linear boundary coditions: $ \frac{\partial u}{\partial\bm{\nu}}=\psi^{-1}(d)$ and $u=\psi^{-1}(d)$, respectively. Thus we can conduct ISS estimates for the considered systems as in \cite{Zheng:2019b} by the sppliting technique combined with the penalty method (see \cite[Remark 5]{Zheng:2019b}).}
     \item {}{The requirement on the smoothness of these functions in \eqref{continuity} and the compatibility condition \eqref{compatibility 1} are only for establishing the existence and regularity of a classical solution of the considered PDEs, and can be weakened for the ISS analysis if weak solutions are considered (see also \cite[Remark 3]{Zheng:2019b}.}

     \item {}{Indeed, we can weaken the condition \eqref{C} to be ``$
         c\geq 0 \ \text{in}\   {B}_R\times \mathbb{R}_{+}
     $''. For example, we consider \eqref{LPE1'} with the Robin boundary condition \eqref{Robin}. Noting that there always exists $\rho\in C^2(\overline{B}_R;\mathbb{R}_+)$ such that
       \begin{align*}
        -\text{div}\ (a\nabla \rho)+\bm{b}\cdot\nabla\rho\geq c_0&>0\ \text{in}\ B_R\times \mathbb{R}_{+},\\
         \frac{\partial\rho}{\partial \bm{\nu}}&\geq 0\ \text{on}\ \partial B_R\times \mathbb{R}_{+},
        \end{align*}
         where  
          $ c_0$ is a positive constant depending on $\rho$.
         Using $u=w\rho$, we can transform the $u$-system \eqref{LPE1'} into $w$-system with the coefficient of $w$, denoted by  $\widehat{c}$, satisfying
  \begin{align*}
         \widehat{c}=c+\frac{1}{\rho}\bm{b}\cdot\nabla\rho-\frac{1}{\rho}\text{div}\ (a\nabla \rho) \geq \frac{c_0}{\max\limits_{\overline{B}_R}\rho}>0.
          \end{align*}
         Moreover, the $w$-system has the structural conditions as \eqref{ABCCBA}, \eqref{H} and \eqref{psi}. Then we can prove the ISS of the $w$-system, which results in the ISS of $u$-system. Due to the spacial limitation, we omit the details.}
     \item  It should be mentioned that proceeding as in this paper and with more specific computations, one may establish ISS estimates for \eqref{LPE1'} over any bounded domain $\Omega \in\mathbb{R}^n$ with a smooth enough boundary.
 \end{enumerate}
 \end{remark}

\section{Maximum Estimate for Parabolic PDEs}\label{Sec: Maximum estimate}
\subsection{Weak maximum principle}
\begin{lemma} \label{weak maximum principle}
Assume that $c\geq 0$ is bounded in $Q_T $ and {}{${}{u\in C^{2,1}(\overline{Q}_T)}$} satisfies $ {}{L_t[u]+N[u]}=f\leq 0$ (resp. $\geq 0$) in $Q_T $, then
\begin{align*}
\max\limits_{\overline{Q}_T}u \leq \max\limits_{\partial_pQ_T}{}{u_+} \ \ \bigg(\text{resp. } \min_{\overline{Q}_T}{}{u_-} \geq \min_{\partial_pQ_T}u\bigg),
\end{align*}
{}{where $u_{+}=\max\{0, u\}$ and $u_{-}=\min\{0, u\}$.}
\end{lemma}
{It seems that the result given in Lemma~\ref{weak maximum principle} is trivial. Nevertheless, for the completeness, we provide a proof by following a similar way given in \cite[page 237]{Wu2006}.}

\begin{pf*}{Proof}
We first show the claim when $f\leq 0$ by contradiction. Indeed, if the claim were false, then there would exist a point $ (x_0,t_0)\in \overline{Q}_T\setminus \partial_pQ_T$ such that
$
u(x_0,t_0)=\max_{\overline{Q}_T} u>0.
$
Thus
\begin{align*}
&u_t|_{(x_0,t_0)}\geq 0,(\nabla u)|_{(x_0,t_0)}=\bm{0},(\text{div} (a\nabla u))|_{(x_0,t_0)}=(a\Delta u)|_{(x_0,t_0)}
\leq 0,c(x_0,t_0)u(x_0,t_0)\geq 0,
\end{align*}
{}{and} $h(x_0,t_0,u(x_0,t_0))>h(x_0,t_0,0)=0
$ due to \eqref{H}.

We have then
\begin{align*}
0\geq f(x_0,t_0)= {}{L_t[u]|_{(x_0,t_0)}}+h(x_0,t_0,u(x_0,t_0))>0,
\end{align*}
which is a contradiction and hence, the claim is valid for $f\leq 0$. For the case of $f\geq 0$, one can proceed in the same way to complete the proof.
\end{pf*}

\subsection{Maximum estimate for parabolic PDEs with a nonlinear Robin boundary condition}
\begin{proposition}\label{maximum estimate}Let ${}{u\in C^{2,1}(\overline{Q}_T)}$ be the solution of the following parabolic equation:
\begin{align*}
 {}{L_t[u]+N[u]}&= {}{f\  \text{in}}\   Q_T,\\
 \frac{\partial u}{\partial\bm{\nu}}+\psi(u)&= {}{d\  \text{on}}\  \partial B_{R} \times (0,T),\\
u {}{(\cdot,0)}&= {}{\phi\  \text{in}}\   B_{R}.
\end{align*}
Then \begin{align*}
\max\limits_{\overline{Q}_T}|u|\leq pR^2+q,
\end{align*}
where
$p= \frac{1}{2 R}\sup\limits_{\partial B_{R} \times (0,T)} |d|$, and $
q= \max\bigg\{\frac{1}{\underline{c}}\bigg( \sup\limits_{Q_T}|f|+2p(\overline{a}n+R\overline{a}+R\overline{b} )\bigg),\sup\limits_{B_{R}}|\phi|\bigg\}.
$
\end{proposition}
\begin{pf*}{Proof}
{}{For $(x,t)\in \overline{Q}_T$, let} $M(x)=p|x|^2+q$ and $v(x,t)=M(x)\pm u(x,t)$. {}{By \eqref{H}, \eqref{A}, \eqref{B} and \eqref{C}, it follows that}
\begin{align*}
 {}{L_t[v]+N[v]}
=&{}{L_t[M]\pm ( L_t[u]+N[u])+N[v]\mp N[u]}\\
=& -2p(an+\nabla a\cdot x)+2p\textbf{b}\cdot x+cp|x|^2+cq{}{\pm f}
{}{+N[M\pm u]\mp N[ u]}\\
\geq &-2p\big(an+R|\nabla a|+R|\textbf{b}| \big)+cq{}{
{\pm f}}
{}{+(N[\pm u]\mp N[ u])}\\
\geq &-2p\big(\overline{a}n+R\overline{a}+R\overline{b} \big)+{}{\underline{c}q
\pm f\ \text{in}\ Q_T}.
\end{align*}
Noting that $
\underline{c}q \geq \sup\limits_{Q_T}|f|+2p\big(\overline{a}n+R\overline{a}+R\overline{b} \big)
$,
we get
\begin{align*}
{}{L_t[v]+N[v]\geq 0\ \text{in}\ Q_T}.
\end{align*}
 By Lemma~\ref{weak maximum principle}, if $v$ has a negative minimum, then $v$ attains the negative minimum on the parabolic boundary $\partial_pQ_T$. On the other hand, noting that $v(x,0)=M(x)\pm \phi (x) \geq 0$ in $B_{R}$, then $v$ attains the negative minimum on $ \partial B_{R} \times (0,T)$, i.e., there exists a point $(x_0,t_0)\in \partial B_{R} \times (0,T)$, such that $v(x_0,t_0) $ is the negative minimum. Thus, $\frac{\partial v}{\partial\bm{\nu}}\big|_{(x_0,t_0)}\leq 0$. Then, at the point $(x_0,t_0)$, {}{we get by \eqref{psi}}
\begin{align*}
0\geq& \frac{\partial v}{\partial\bm{\nu}}+\psi(0)\\
>& \frac{\partial v}{\partial\bm{\nu}}+\psi(v)\\
= &\frac{\partial M}{\partial\bm{\nu}}+\psi(M\pm u)\mp\psi(u)\pm\bigg(\frac{\partial u}{\partial\bm{\nu}}+\psi(u)\bigg)\\
= &
2p{R}+(\psi(M\pm u){}{\mp\psi( u)})\pm d\\
\geq &{}{2p{R}+\big(\psi(\pm u)\mp\psi( u)}\big) \pm d\\
\geq &{2p}{R} \pm d\\
=& \sup\limits_{\partial B_{R} \times (0,T)} |d| \pm d\\
\geq& 0,
\end{align*}
which is a contradiction. Therefore, there must be $v\geq0$ in $\overline{Q}_T $, which follows {}{that $
|u|\leq M\leq pR^2+q
$ in $\overline{Q}_T$.}
\end{pf*}

\subsection{Maximum estimate for parabolic PDEs with a nonlinear Dirichlet boundary condition}

\begin{proposition}\label{Dirichlet'}
Let ${}{u\in C^{2,1}(\overline{Q}_T)}$ be the solution of the following parabolic equation:
\begin{subequations}\label{LPE2601}
\begin{align}
{}{L_t[u]+N[u]}&={}{f\  \text{in}}\   Q_T,\\
 \psi(u)&={}{d\  \text{on}}\  \partial B_{R} \times (0,T),\\
{}{u(\cdot,0)}&={}{\phi\  \text{in}}\   B_{R}.
\end{align}
\end{subequations}
Then
\begin{align*}
\max\limits_{\overline{Q}_T}|u|\leq \max\bigg\{\frac{1}{\underline{c}}\sup\limits_{Q_T}|f|,{\psi}^{-1}\bigg(\sup\limits_{\partial B_{R} \times (0,T)}|d|\bigg),\sup\limits_{B_{R}}|\phi|\bigg\}.
\end{align*}
\end{proposition}
\begin{pf*}For $(x,t)\in \overline{Q}_T$, let  $v(x,t)=M\pm u(x,t)$ with $M=\max\bigg\{\frac{1}{\underline{c}} \sup\limits_{Q_T}|f|,{\psi}^{-1}\bigg(\sup\limits_{\partial B_{R} \times (0,T)}|d|\bigg),\sup\limits_{B_{R}}|\phi|\bigg\}$.

Proceeding as in the proof of Proposition~\ref{maximum estimate}, we have
\begin{align*}
{}{L_t[v]+N[v]\geq 0\ \text{in}\ Q_T}.
\end{align*}
Then it suffices to show that if $v(x_0,t_0) $ is the negative minimum at some point $(x_0,t_0)\in \partial B_{R} \times (0,T)$, we will obtain a contradiction. Indeed, noting that
\begin{align*}
 \psi(u(x_0,t_0))=d(x_0,t_0) \leq \sup\limits_{\partial B_{R} \times (0,T)}|d|,
\end{align*}
it follows that
\begin{align*}
u(x_0,t_0)\leq  {\psi}^{-1}\bigg(\sup\limits_{\partial B_{R} \times (0,T)}|d|\bigg)\leq M.
\end{align*}
Then we have
$
0>v(x_0,t_0) =M\pm u(x_0,t_0)\geq 0,
$
which is actually a contradiction.
\end{pf*}

\section{EISS Estimates for Parabolic PDEs with Different Types of Nonlinear Boundary Conditions}\label{Sec: EISS estimates}
\begin{pf*}{Proof of Theorem \ref{main result}} We proceed on the proof in the following 3~steps.

(i) We establish an EISS estimate of the {}{solution} to \eqref{LPE1'} with the nonlinear Robin boundary condition~\eqref{Robin}.

Let $v\in {}{C^{2,1}(\overline{Q}_T)}$ be the unique solution of the following parabolic equation:
\begin{align*}
{}{L_t[v]+N[v]}&={}{f\  \text{in}}\   Q_T,\\
 \frac{\partial v}{\partial\bm{\nu}}+\psi(v)&={}{d\  \text{on}}\  \partial B_{R} \times (0,T),\\
{}{v(\cdot,0)}&=0{}{\  \text{in}}\   B_{R}.
%
\end{align*}

According to Proposition~\ref{maximum estimate}, we have
\begin{align*}
\max\limits_{\overline{Q}_T}|v|\leq{}{R_0} \sup\limits_{\partial B_{R} \times (0,T)} |d|+\frac{1}{\underline{c}} \sup\limits_{Q_T}|f|,
\end{align*}
{}{where $R_0=\frac{R}{2}+ \frac{1}{\underline{c}R}( \overline{a}n+R\overline{a}+R\overline{b} )$.}

Let $w=u-v$. It is obvious that $w$ satisfies:
\begin{subequations}\label{LPE2}
\begin{align}
{}{L_t[w]+N[u]-N[v]}&=0\  {}{\text{in}}\   Q_T,\\
 \frac{\partial w}{\partial\bm{\nu}}+\psi(u)-\psi(v)&=0\  {}{\text{on}}\    \partial B_{R} \times (0,T),\\
{}{w(\cdot,0)}&={}{\phi\  \text{in}}\   B_{R}.
\end{align}
\end{subequations}
Multiplying \eqref{LPE2} with $w$ and integrating by parts, we have
\begin{align*}
&\frac{1}{2}\frac{\text{d}}{\text{d}t}\|w\|^2+\|\sqrt{a}\nabla w\|^2\\
=&-\int_{B_R}cw^2\text{d}x-\int_{B_R} (\textbf{b}  \cdot \nabla w)w\text{d}x
+\int_{\partial B_R}aw\nabla w\cdot \bm{\nu}\text{d}S+\int_{B_R}{}{(N[v](x,t)-N[u](x,t))}(u-v)\text{d}x.
\end{align*}
Applying the formula of integration by parts, {}{the Trace~Theorem (see the appendix) and by \eqref{B},} we have
\begin{align*}
-\int_{B_R} (\textbf{b}  \cdot \nabla w)w\text{d}x
=&\frac{1}{2}\int_{B_R}(\text{div} \ \bm{b})w^2\text{d}x- \frac{1}{2}\int_{\partial B_R}w^2( \textbf{b}\cdot \bm{\nu})\text{d}S\\
\leq&\frac{\overline{b}}{2}\left(\int_{B_R}w^2\text{d}x+ \int_{\partial B_R}w^2\text{d}S\right)
\\
\leq &\frac{\overline{b}}{2}\left(\|w\|^2+ C^2_{\text{Trace}}(\|w\|+\|\nabla w\|)^2\right)\\
\leq &\frac{\overline{b}}{2}\left( (1+2C^2_{\text{Trace}})\|w\|^2+2C^2_{\text{Trace}} \|\nabla w\|^2\right).
\end{align*}
{}{By \eqref{psi} and \eqref{H},} we always have
\begin{align*}
&\int_{\partial B_R}aw\nabla w\cdot \bm{\nu}\text{d}S=\int_{\partial B_R}a(\psi(v)-\psi(u))(u-v)\text{d}S\leq 0,\\
&\int_{B_R}{}{(N[v](x,t)-N[u](x,t))}(u-v)\text{d}x\leq 0.
\end{align*}
Thus, we obtain {}{by \eqref{A}, \eqref{C} and \eqref{ABC}}
\begin{align*}
\frac{d}{dt}\|w\|^2\leq& -2\|\sqrt{a}\nabla w\|^2-2\|\sqrt{c} w\|^2
+\overline{b}\left( (1+2C^2_{\text{Trace}})\|w\|^2+2C^2_{\text{Trace}} \|\nabla w\|^2\right)\\
\leq&
-\left(2\underline{c}-\overline{b}(1+2C^2_{\text{Trace}})\right)\| w\|^2
+ 2 (\overline{b}C^2_{\text{Trace}}-\underline{a}) \|\nabla w\|^2\\
\leq &-\left(2\underline{c}-\overline{b}(1+2C^2_{\text{Trace}})\right)\| w\|^2,
\end{align*}
which gives
$
\|w(\cdot,T)\|^2\leq\|\phi\|^2\e^{ -\big(2\underline{c}-\overline{b}\left(1+2C^2_{\text{Trace}}\right)\big)T}.
$

Finally, we have
\begin{align}
\|u(\cdot,T)\|
\leq &\|w(\cdot,T)\|+\|v(\cdot,T)\|\notag\\
\leq& \|\phi\|\e^{ -\big(\underline{c}-2\overline{b}\left(1+2C^2_{\text{Trace}}\right)\big)T}
+ {}{R_0}\sqrt{|B_R|}\sup\limits_{\partial B_{R} \times (0,T)} |d|
+\frac{1}{\underline{c}} \sqrt{|B_R|}\sup\limits_{Q_T}|f|,\ \ \forall T>0.\label{2801}
\end{align}

(ii) We establish an EISS estimate of the {}{solution} to \eqref{LPE1'} with the nonlinear Neumann boundary condition~\eqref{Neumann}.
Let $v\in {}{C^{2,1}(\overline{Q}_T)}$ be the unique solution of the following parabolic equation:
\begin{align*}
{}{L_t[v]+N[v]}&={}{f\  \text{in}}\   Q_T,\\
 \frac{\partial v}{\partial\bm{\nu}}+v&=\psi^{-1}(d)\  {}{\text{on}}\   \partial B_{R} \times (0,T),\\
{}{v(\cdot,0)}&=0\  {}{\text{in}}\  B_{R}.
\end{align*}

According to Proposition~\ref{maximum estimate}, we have
\begin{align}
\max\limits_{\overline{Q}_T}|v|\leq&  {}{R_0}\sup\limits_{\partial B_{R} \times (0,T)} |\psi^{-1}(d)|+\frac{1}{\underline{c}} \sup\limits_{Q_T}|f|\notag\\
\leq& {}{R_0}\psi^{-1}\bigg(\sup\limits_{\partial B_{R} \times (0,T)} |d|\bigg) +\frac{1}{\underline{c}} \sup\limits_{Q_T}|f|,\label{022703}
\end{align}
{}{where $R_0=\frac{R}{2}+ \frac{1}{\underline{c}R}( \overline{a}n+R\overline{a}+R\overline{b} )$}.

Let $w=u-v$, which satisfies:
\begin{align*}
{}{L_t[w]+N[u]-N[v]}&=0\  {}{\text{in}}\  Q_T,\\
 \frac{\partial w}{\partial\bm{\nu}}-v&=0\  {}{\text{on}}\   \partial B_{R} \times (0,T),\\
{}{w(\cdot,0)}&={}{\phi\  \text{in}}\  B_{R}.
\end{align*}
{}{By Young's inequality, the Trace Theorem and \eqref{A}, it follows that}
\begin{align*}
 \int_{\partial B_R}aw\nabla w\cdot \bm{\nu}\text{d}S 
=&\int_{\partial B_R}awv\text{d}S\\
\leq & \frac{\varepsilon }{2}\int_{\partial B_R}w^2\text{d}S +\frac{1}{2\varepsilon}\int_{\partial B_R}a^2v^2\text{d}S,\\
\leq & \varepsilon C^2_{\text{Trace}}\left( \|w\|^2+\|\nabla w\|^2\right) + \frac{1}{2\varepsilon} \overline{a}^2  |\partial B_R|\max\limits_{{}{\partial B_R}}|v|^2.
\end{align*}

Proceeding in the same way as in (i), we get
\begin{align*}
\frac{\text{d}}{\text{d}t}\|w\|^2\leq&
-\left(2\underline{c}-\overline{b}(1+2C^2_{\text{Trace}})-2\varepsilon C^2_{\text{Trace}}\right)\| w\|^2
+ 2 \left(\overline{b}C^2_{\text{Trace}}-\underline{a}-\varepsilon C^2_{\text{Trace}}\right) \|\nabla w\|^2
+\frac{1}{\varepsilon} \overline{a}^2  |\partial B_R|\max\limits_{{}{{}{\partial B_R}}}|v|^2\notag\\
\leq &-\left(2\underline{c}-\overline{b}(1+2C^2_{\text{Trace}})-2\varepsilon C^2_{\text{Trace}}\right)\| w\|^2
+\frac{1}{\varepsilon} \overline{a}^2  |\partial B_R|\max\limits_{ \overline{B}_R}|v|^2\notag\\
:=&-\lambda \|w\|^2+V(t),
\end{align*}
where we choose $ \varepsilon>0$ small enough such that $ \lambda=2\underline{c}-\overline{b}(1+2C^2_{\text{Trace}})-2\varepsilon C^2_{\text{Trace}}>0$ and $\overline{b}C^2_{\text{Trace}}-\underline{a}-\varepsilon C^2_{\text{Trace}}>0$ {}{due to \eqref{ABC}}.

By Gronwall's inequality, we have
\begin{align}
\|w(\cdot, T)\|^2\leq& \|\phi\|^2\e^{-\lambda T}+\max\limits_{t\in [0,T]}V(t)\cdot \int_{0}^T \e^{-\lambda (T-t)}\text{d}t\notag\\
\leq &
\|\phi\|^2\e^{-\lambda T}+\frac{1}{\lambda}\max\limits_{t\in [0,T]}V(t)\notag\\
\leq &\|\phi\|^2\e^{-\lambda T} +\frac{\overline{a}^2}{\varepsilon\lambda}   |\partial B_R|\max\limits_{{\overline{Q}_T}}|v|^2.\label{022704}
\end{align}
Finally, by $\|u(\cdot,T)\|
\leq \|w(\cdot,T)\|+\|v(\cdot,T)\|$, \eqref{022703} and \eqref{022704}, for any $T>0$, it follows {}{that}
\begin{align}
\|u(\cdot, T)\|\leq&
\|\phi\|\e^{-\frac{\lambda}{2} T}+{}{R_0}\bigg(1+ \frac{\overline{a}}{\sqrt{\varepsilon\lambda}}\bigg)\sqrt{|\partial B_R|}\psi^{-1}\bigg(\sup\limits_{\partial B_{R} \times (0,T)} |d|\bigg)+\frac{1}{\underline{c}}\bigg(1+ \frac{\overline{a}}{\sqrt{\varepsilon\lambda}}\bigg)\sqrt{|\partial B_R|} \sup\limits_{Q_T}|f|.
\label{2802}
\end{align}

(iii) For the EISS estimate of the {}{solution} to \eqref{LPE1'} with the nonlinear Dirichlet boundary condition~\eqref{Dirichlet}, it suffices to estimate the solutions of the following parabolic equations:
\begin{align*}
{}{L_t[v]+N[v]}&={}{f\  \text{in}}\   Q_T,\\
 v&=\psi^{-1}(d)\  {}{\text{on}}\   \partial B_{R} \times (0,T),\\
{}{v(\cdot,0)}&=0\  {}{\text{in}}\   B_{R},
\end{align*}
and
\begin{align*}
{}{L_t[w]+N[u]-N[v]}&=0\  {}{\text{in}}\    Q_T,\\
 w&=0\  {}{\text{on}}\  \partial B_{R} \times (0,T),\\
{}{w(\cdot,0)}&={}{\phi\  \text{in}}\   B_{R},
\end{align*}
where $w=u-v$.

Indeed, by Proposition~\ref{Dirichlet'}, we have
\begin{align*}
\max\limits_{\overline{Q}_T}|v|\leq \max\bigg\{\frac{1}{\underline{c}}\sup\limits_{Q_T}|f|,{\psi}^{-1}\bigg(\sup\limits_{\partial B_{R} \times (0,T)}|d|\bigg)\bigg\}.
\end{align*}
Proceeding as in (i), we get
\begin{align*}
\|w(\cdot,T)\|\leq \e^{-\frac{1}{2}\left(2\underline{c}-\overline{b}\right)T} \|\phi\|,\ \ \forall T>0.
\end{align*}
Finally, for any $T>0$, it {}{follows that}
\begin{align}
\|u(\cdot, T)\|
\leq  \|\phi\|\e^{-\frac{1}{2}\left(2\underline{c}-\overline{b}\right)T}
+\sqrt{|B_{R}|}\max\bigg\{\frac{1}{\underline{c}}\sup\limits_{Q_T}|f|,{\psi}^{-1}\bigg(\sup\limits_{\partial B_{R} \times (0,T)}|d|\bigg)\bigg\}.
\label{2803}
\end{align}
\end{pf*}
\section{An Illustrative Example}\label{Sec: example}
{}{ We consider the following super-linear parabolic equation:
\begin{align}\label{Example}
u_t-\text{div} \ (a\nabla u)+cu+u\ln(1+u^2)=f\ \text{in}\  B_R\times \mathbb{R}_+,\
\end{align}
coupled with the nonlinear Robin boundary condition:
\begin{align}\label{Example 1a}
\frac{\partial u}{\partial\bm{\nu}}+u+u^3=d\ \ \text{on} \  \partial B_{R} \times \mathbb{R}_+,
\end{align}
or the nonlinear Dirichlet boundary condition:
\begin{align}
u+u^3=d\ \ \text{on} \  \partial B_{R} \times \mathbb{R}_+. \label{Example 1b}
\end{align}
The initial value condition is given by:
\begin{align*}
  u{}{(\cdot,0)=\phi(\cdot)}\ \ \text{in} \  B_R.
\end{align*}
{}{Note that $N[u](x,t)=h(x,t,u)\equiv u\ln(1+u^2)$ and $\psi(u)=u+u^3$, both of which are $C^2$-continuous, odd and strictly increasing in $u$. Thus \eqref{H} and \eqref{psi} are satisfied.}
%

{}{If we assume that $a$, $c$, $f$, $d$, and $\phi$ satisfy \eqref{continuity}, \eqref{ABCCBA} and \eqref{compatibility 1}, then according to Theorem~\ref{main result}, system~\eqref{Example} with \eqref{Example 1a}, or \eqref{Example 1b}, is EISS, having the following estimate for any $T>0$:
\begin{align*}
\|u(\cdot,T)\|
\leq \|\phi\|\e^{ -\underline{c}T}+ \sqrt{|B_R|}G\bigg(\sup\limits_{Q_T}|f|,\sup\limits_{\partial B_{R} \times (0,T)} |d|\bigg),
\end{align*}
where $G(y,z)=\frac{1}{\underline{c}}y+\left(\frac{R}{2}+ \frac{\overline{a}}{\underline{c}} (n+R) \right)z,$
 or $G(y,z)=\max\left\{\frac{1}{\underline{c}}y,{\psi}^{-1}(z)\right\}.
 $}
\section{Concluding Remarks}\label{Sec: Conclusion}
This paper presented an application of the maximum principle-based approach proposed in \cite{Zheng:2019a,Zheng:2019b} to the establishment of ISS properties w.r.t. in-domain and boundary disturbances for certain nonlinear parabolic PDEs over higher dimensional domains with different types of nonlinear boundary conditions. The proposed scheme for achieving the ISS estimates of the solution is based on the the Lyapunov method and the maximum estimates for parabolic PDEs with nonlinear boundary conditions. An ISS analysis for a parabolic PDE with a super-linear term and nonlinear boundary conditions has been carried out, which demonstrated the {}{effectiveness} of the developed approach.


\section{Appendix: Trace Theorem}
 Let $H^1(B_R):=\{ u: B_R\rightarrow \mathbb{R}\ \text{is\ locally\ summable}| u \in L^2(B_R),\nabla u\in( L^2(B_R))^n\}$ endowed with the norm $\|u\|_{H^1(B_R)}:= \|u\|_{L^2(B_R)}+\|\nabla u\|_{L^2(B_R)}$.
\begin{theorem} (Trace Theorem \cite[Chapter 5]{Evans:2010}) There exists a bounded linear operator $\mathcal {T}: H^1(B_R)\rightarrow L^2(\partial B_R)$ such that
\begin{itemize}
\item[(i)] $\mathcal {T}u=u|_{\partial B_R}$ if $ u\in H^1(B_R)\cap C(\overline{B}_R)$,
and \item[(ii)] $ \|\mathcal {T}u\|_{L^2(\partial B_R)}\leq C_{\text{Trace}}\|u\|_{H^1(B_R)}$ for each $u\in H^1(B_R)$, with the constant $  C_{\text{Trace}}$ depending only on $B_R$.
\end{itemize}
\end{theorem}

\end{document}